\documentclass[prodmode,acmtoms]{acmsmall} 
\usepackage{latexsym,amsmath,amssymb,xspace,comment,graphicx, url, color}
\usepackage[ruled]{algorithm2e}
\renewcommand{\algorithmcfname}{ALGORITHM}
\SetAlFnt{\small}
\SetAlCapFnt{\small}
\SetAlCapNameFnt{\small}
\SetAlCapHSkip{0pt}
\IncMargin{-\parindent}

\newlength\bshft
\def\fakebold#1{\ooalign{$#1$\cr\kern-\bshft$#1$\cr\kern\bshft$#1$}}

\begin{document}

\markboth{S. Harase and T. Kimoto}{Implementing 64-bit Maximally Equidistributed $\mathbb{F}_2$-Linear Generators}

\title{Implementing 64-bit Maximally Equidistributed $\mathbb{F}_2$-Linear Generators with Mersenne Prime Period}
\author{SHIN HARASE
\affil{Ritsumeikan University}
TAKAMITSU KIMOTO
\affil{Recruit Holdings Co., Ltd.}}

\begin{abstract}
CPUs and operating systems are moving from 32 to 64 bits, 
and hence it is important to have good pseudorandom number generators designed to fully exploit these word lengths. 
However, existing 64-bit very long period generators based on linear recurrences modulo 2 are not completely optimized in terms of the equidistribution properties. 
Here we develop 64-bit maximally equidistributed pseudorandom number generators 
that are optimal in this respect and have speeds equivalent to 64-bit Mersenne Twisters. 
We provide a table of specific parameters with period lengths from $2^{607}-1$ to $2^{44497}-1$.
\end{abstract}



\acmformat{Shin Harase and Takamitsu Kimoto. 2017. Implementing 64-bit maximally equidistributed 
$\mathbb{F}_2$-linear generators with Mersenne prime period.}

\begin{bottomstuff}
This work was partially supported by the Grants-in-Aid for JSPS Fellows 24$\cdot$7985, for Young Scientists (B) $\#$26730015, for Scientific Research (B) $\#$26310211, and for Challenging Exploratory Research $\#$15K13460 from the Japan Society for the Promotion of Scientific Research. This work was also supported by  JST CREST.
Authors' addresses: S. Harase,
College of Science and Engineering, Ritsumeikan University, 1-1-1 Nojihigashi, Kusatsu, Shiga, 525-8577, Japan; e-mail: harase@fc.ritsumei.ac.jp; 
T. Kimoto,  Recruit Holdings Co., Ltd., 1-9-2 Marunouchi, Chiyoda-ku, Tokyo, 100-6640, Japan; tkimoto@r.recruit.co.jp.
\end{bottomstuff}

\maketitle

\section{Introduction}
Monte Carlo simulations are a basic tool in financial engineering, computational physics, statistics, and other fields. 
To obtain precise simulation results, the quality of pseudorandom number generators is important. 
At present, the 32-bit Mersenne Twister (MT) generator MT19937 (with period $2^{19937}-1$) \cite{MT19937} is one of the most widely used pseudorandom number generators. 
However, modern CPUs and operating systems are moving from 32 to 64 bits, 
and hence it is important to have high-quality generators designed to fully exploit 64-bit words. 

Many pseudorandom number generators, including Mersenne Twisters, are based on linear recurrences modulo 2; 
these are called {\it $\mathbb{F}_2$-linear generators}.  
One advantage of these generators is that they can be assessed by means of the {\it dimension of equidistribution with $v$-bit accuracy}, 
which is a most informative criterion for high dimensional uniformity of the output sequences. 
In fact, MT19937 is not completely optimized in this respect. 
\citeN{PannetonLM06} developed the Well Equidistributed Long-period Linear (WELL) generators with periods from $2^{512}-1$ to $2^{44497}-1$, 
which are completely optimized for this criterion (called {\it maximally equidistributed}), 
but the parameter sets were only searched for the case of 32-bit generators. 
Conversely, there exist several 64-bit $\mathbb{F}_2$-linear generators. 
\citeN{Nishimura2000} developed 64-bit Mersenne Twisters, 
and the SIMD-oriented Fast Mersenne Twister (SFMT) generator \cite{MR2743921} has a function to generate 64-bit unsigned integers. 
For graphics processing units, Mersenne Twister for Graphic Processors (MTGP) \cite{MR3031631} is also a good candidate. 
However, these generators are not maximally equidistributed. 
In earlier work, \citeN{MR1620231} searched for 64-bit maximally equidistributed combined Tausworthe generators (with some additional properties). 
At present, though, to the best of our knowledge, there exists no 64-bit maximally equidistributed MT-type $\mathbb{F}_2$-linear generator with period $2^{19937}-1$, such as a 64-bit variant of the WELL generators.

The aim of this article is to develop 64-bit maximally equidistributed $\mathbb{F}_2$-linear generators 
with similar speed as the 64-bit Mersenne Twisters \cite{Nishimura2000}. 
The key techniques are (i) state transitions with double feedbacks \cite{PannetonLM06,MR2743921} 
and (ii) linear output transformations with several memory references \cite{Harase2009}. 
We provide a table of specific parameters with periods from $2^{607}-1$ to $2^{44497}-1$. 
The design of our generators is based on a combination of existing techniques, such as the WELL and dSFMT generators \cite{PannetonLM06,MR2743921}, 
but we select state transitions differently from those of the original WELL generators to maintain the generation speed. 
We refer to these generators as {\it 64-bit Maximally Equidistributed $\mathbb{F}_2$-Linear Generators (MELGs) with Mersenne prime period}.

In practice, we often convert unsigned integers into 53-bit double-precision
floating-point numbers in $[0, 1)$ in IEEE 754 format. Our 64-bit generators are useful for this. 
To generate 64-bit output values, one can either use a pseudorandom number generator whose linear recurrence is implemented with 32-bit integers 
and then take two successive 32-bit blocks or, instead, use a pseudorandom number generator whose recurrence is implemented directly over 64-bit integers. 
As described below, the former method may degrade the dimension of equidistribution with $v$-bit accuracy, 
compared with simply using 32-bit output values. We consider the case of the 32-bit MT19937 generator in Section~\ref{sec:performance}.  
For this reason, we develop 64-bit MELGs to directly generate 64-bit unsigned integers.

The article is organized as follows. In the next section, 
we review $\mathbb{F}_2$-linear generators and their theoretical criteria. 
We also summarize the framework of Mersenne Twisters. 
Section~\ref{sec:MELG} is devoted to our main result: 64-bit MELGs. 
In Section~\ref{sec:performance}, we compare our generators with others in terms of 
speeds, theoretical criteria, and empirical statistical tests. 
Section~\ref{sec:conclusions} concludes. 

\section{Preliminaries} \label{sec:preliminaries}
\subsection{$\protect\fakebold{\mathbb{F}}_{\mathbf{2}}$-Linear Generators} \label{subsec:F2-linear generators}

We recall the notation of $\mathbb{F}_2$-linear generators; see \cite{LP2009,MatsumotoSHN06} for details. 
Let $\mathbb{F}_2:= \{ 0, 1\}$ be the two-element field, i.e., 
addition and multiplication are performed modulo 2. We consider the following class of generators.
\begin{definition}[$\mathbb{F}_2$-linear generator] \label{def:F2-linear}
Let $S := \mathbb{F}_2^p$ be a $p$-dimensional state space (of the possible states of the memory assigned for generators).  
Let $f : S \to S$ be an $\mathbb{F}_2$-linear state transition function.  
Let $O := \mathbb{F}_2^w$ be the set of outputs, where $w$ is the word size of the intended machine, 
and let $o: S \to O$ an $\mathbb{F}_2$-linear output function. 
For an initial state $s_0 \in S$, at every time step, the state is changed by the recursion 
\begin{eqnarray} \label{eqn:transition}
s_{i+1} = f(s_{i}) \quad (i = 0, 1, 2, \ldots ),
\end{eqnarray} 
and the output sequence is given by 
\begin{eqnarray} \label{eqn:output}
o(s_0), o(s_1), o(s_2), \ldots \in O.
\end{eqnarray}
We identify $O$ as a set of unsigned $w$-bit binary integers. 
A generator with these properties is called an {\it $\mathbb{F}_2$-linear generator}. 
\end{definition}

Let $P(z)$ be the characteristic polynomial of $f$. 
The recurrence (\ref{eqn:transition}) has a period of length $2^p-1$ 
(its maximal possible value) if and only if $P(z)$ is a primitive polynomial modulo 2 \cite{Niederreiter:book,Knuth:1997:ACP:270146}. 
When this value is reached, 
we say that the $\mathbb{F}_2$-linear generator has {\it maximal period}. 
Unless otherwise noted, we assume throughout that this condition holds. 

\subsection{Quality Criteria} \label{subsec:criteria}

Following \cite{LP2009,MatsumotoSHN06}, we recall two quality criteria for $\mathbb{F}_2$-linear generators. 
A most informative criterion for high dimensional uniformity is the {\it dimension of equidistribution with $v$-bit accuracy}.  
Assume that an $\mathbb{F}_2$-linear generator has the maximal period $2^p-1$.
We identify the output set $O:=\mathbb{F}_2^w$ as a set of unsigned $w$-bit binary integers. 
We focus on the $v$ most significant bits of the output, 
and regard these bits as the {\em output with $v$-bit accuracy}.
This amounts to considering the composition 
$o_v: S \stackrel{o}{\to} \mathbb{F}_2^w \to \mathbb{F}_2^v$,
where the latter mapping denotes taking the $v$ most significant bits.
We define the $k$-tuple output function as 
\[
 o_v^{(k)}:S \to (\mathbb{F}_2^v)^k, \quad s_0 \mapsto 
 (o_v(s_0), o_v(f(s_0)), \ldots, o_v(f^{k-1}(s_0))).
\]
Thus, $o_v^{(k)}(s_0)$ is the vector formed by the $v$ most significant bits of $k$ consecutive output values of the pseudorandom number generators from a state $s_0$.
\begin{definition}[Dimension of equidistribution with $v$-bit accuracy] \label{def:equidistribution}
The generator is said to be $k$-dimensionally equidistributed
with $v$-bit accuracy if and only if $o_v^{(k)}:S \to (\mathbb{F}_2^v)^k$ is surjective. 
The largest value of $k$ with this property is called the {\it dimension of equidistribution 
with $v$-bit accuracy}, denoted by $k(v)$.
\end{definition} 
Because $o_v^{(k)}$ is $\mathbb{F}_2$-linear, $k$-dimensional equidistribution with $v$-bit accuracy 
means that every element in $(\mathbb{F}_2^v)^k$ occurs with the same probability,
when the initial state $s_0$ is uniformly distributed over the state space $S$. 
As a criterion of uniformity, larger values of $k(v)$ for each $1 \leq v \leq w$ are desirable \cite{Tootill}. 
We have a trivial upper bound $k(v) \leq \lfloor p/v \rfloor$. 
The gap $d(v) := \lfloor p/v \rfloor -k(v)$ is called the {\it dimension defect at $v$}, and the sum of the gaps 
$\Delta:=\sum_{v = 1}^{w} d(v)$ is called the {\it total dimension defect}. 
If $\Delta = 0$, the generator is said to be {\it maximally equidistributed}. 
For $\mathbb{F}_2$-linear generators, 
one can quickly compute $k(v)$ for $v = 1, \ldots, w$ by using lattice reduction algorithms over formal power series fields \cite{HMS2010,Shin2011141}. 
These are closely related to the lattice reduction algorithm originally proposed by \citeN{CL2000}. 

As another criterion, the number $N_1$ of nonzero coefficients 
for $P(z)$ should be close to $p/2$ \cite{springerlink:10.1007/BF01029989,MR1155576}. 
For example, generators for which $P(z)$ is a trinomial or pentanomial fail in statistical tests 
\cite{Lindholm1968,Matsumoto:1996:SDR:232807.232815,MR1958868}. 
If $N_1$ is not large enough, the generator suffers a long-lasting  impact 
for poor initialization known as a 0-excess state $s_0 \in S$, which contains only a few bits set to 1 \cite{PannetonLM06}.
Thus, $N_1$ should be in the vicinity of $p/2$. 

\subsection{Mersenne Twisters} \label{subsec:MT}

In general, to ensure a maximal period, testing the primitivity of $P(z)$ is the bottleneck in searching for long-period generators 
(because the factorization of $2^p-1$ is required). 
If $p$ is a Mersenne exponent (i.e., $2^p-1$ is a Mersenne prime), 
one can instead use an irreducibility test that is easier and equivalent. 
\citeN{MT19937} developed a pseudorandom number generator with a Mersenne prime period $2^p-1$ known as the Mersenne Twister. 
We briefly review their method. 
The state space $S$ and state transition $f: S \to S$ are expressed as 
\begin{eqnarray} \label{eqn:MT_state_transition} 
(\overline{\mathbf{w}}_{i}^{w-r}, \mathbf{w}_{i+1}, \mathbf{w}_{i+2}, \ldots, \mathbf{w}_{i+N-1}) \mapsto  (\overline{\mathbf{w}}_{i+1}^{w-r}, \mathbf{w}_{i+2}, \mathbf{w}_{i+3}, \ldots, \mathbf{w}_{i+N}),
\end{eqnarray}
where $\mathbf{w}_i \in \mathbb{F}^w$ is a $w$-bit word vector, 
$\overline{\mathbf{w}}_i^{w-r} \in \mathbb{F}_2^{w-r}$ denotes the $w-r$ most significant bits of $\mathbf{w}_i$, 
$N:= \lceil p/w \rceil$, and  $r$ is a non-negative integer such that $p = Nw-r$, 
so that the $p$ bits of $S$ are stored in an array of $Nw$ bits in which there are $r$ unused bits. 
The state transition (\ref{eqn:MT_state_transition}) is implemented as 
\begin{eqnarray} \label{eqn:MT_recurtion2}
\mathbf{w}_{i+N} := \mathbf{w}_{i+M} \oplus (\overline{\mathbf{w}}_i^{w-r} \mid \underline{\mathbf{w}}_{i+1}^{r}) A,
\end{eqnarray}
where $\underline{\mathbf{w}}_{i+1}^{r} \in \mathbb{F}_2^r$ represents the $r$ least significant bits of $\mathbf{w}_{i+1}$, 
$\oplus$ is the bitwise exclusive-OR (i.e., addition in $\mathbb{F}_2^w$), 
and $(\overline{\mathbf{w}}_i^{w-r} \mid \underline{\mathbf{w}}_{i+1}^{r}) $, a $w$-bit vector, is the concatenation of  
the $(w-r)$-bit vector $\overline{\mathbf{w}}_i^{w-r}$ and the $r$-bit vector $\underline{\mathbf{w}}_{i+1}^{r}$ in that order. 
$M$ is an integer such that $0 < M < N-1$, and $A \in \mathbb{F}_2^{w \times w}$ is a $(w \times w)$-regular matrix 
(with the format in (\ref{eqn:def_mat_A}) below). 
Furthermore, to improve $k(v)$, for the right-hand side of (\ref{eqn:MT_state_transition}), 
the output transformation $o: S \to O$ is implemented as 
\begin{eqnarray} \label{eqn:tempering}
(\overline{\mathbf{w}}_{i+1}^{w-r}, \mathbf{w}_{i+2}, \ldots, \mathbf{w}_{i+N}) \in S \mapsto \mathbf{w}_{i + N}T \in O, 
\end{eqnarray}
where $T \in \mathbb{F}_2^{w \times w}$ is a suitable $(w \times w)$-regular matrix. 
This technique is called {\it tempering}. From this, we obtain an output sequence $\mathbf{w}_{N}T, \mathbf{w}_{N+1}T, \mathbf{w}_{N+2}T, \ldots \in O$ by multiplying a matrix $T$ by the sequence from (\ref{eqn:MT_recurtion2}). 
\citeN{MT19937} and \citeN{Nishimura2000} searched for 
32- and 64-bit Mersenne Twisters with period length $2^{19937}-1$, respectively. 
In Appendix~B of \cite{MT19937}, it is proved that these generators cannot attain 
the maximal equidistribution (e.g., $\Delta = 6750$ for MT19937 in \cite{MT19937}). 
In fact, the state transition in (\ref{eqn:MT_state_transition})--(\ref{eqn:MT_recurtion2}) is very simple,  
but the linear output transformation $T$ in (\ref{eqn:tempering}) is rather complicated (see \cite{MT19937,Nishimura2000} for details).

\section{Main Result: 64-bit Maximally Equidistributed $\mathbb{F}_2$-Linear Generators} \label{sec:MELG}
\subsection{Design} \label{sibsec:design}

To obtain maximally equidistributed generators without loss of speed, we try to shift the balance of costs in $f$ and $o$. 
The key technique is a suitable choice of (i) state transitions with double feedbacks proposed in \cite{PannetonLM06,MR2743921} and 
(ii) linear output transformations with several memory references from \cite{Harase2009}. 
Let $N = \lceil p/w \rceil$. We divide $S = \mathbb{F}_2^p$ into two parts $S = \mathbb{F}_2^{p-w} \times \mathbb{F}_2^w$ 
and consider the state transition $f: S \to S$ with 
\begin{eqnarray} \label{eqn:MELG_state_transition}
(\overline{\mathbf{w}}_i^{w-r}, \mathbf{w}_{i+1}, \mathbf{w}_{i+2}, \ldots, \mathbf{w}_{i+N-2}, \mathbf{v}_i) \mapsto (\overline{\mathbf{w}}_{i+1}^{w-r}, \mathbf{w}_{i+2}, \mathbf{w}_{i+3}, \ldots, \mathbf{w}_{i+N-1}, \mathbf{v}_{i+1}), 
\end{eqnarray}
where $\mathbf{w}_{i+N-1}$ and $\mathbf{v}_{i+1}$ are determined by the recursions (\ref{eqn:MELG_recurtion1_def}) and (\ref{eqn:MELG_recurtion2_def}), 
as described below. The first $p-w$ bits are stored in an array in which $r$ bits are unused, as for the original MTs. 
Note that the number of words is $N-1$, not $N$. 
The remaining word $\mathbf{v}_i$ is expected to be stored in a register of the CPU and updated as $\mathbf{v}_{i+1}$ at the next step, 
so that the implementation requires only a single word 
(see Figure~1 and Algorithm~1 in this subsection). 
The use of the extra state variable $\mathbf{v}_i$ was originally proposed by \citeN{PannetonLM06} and refined by \citeN{MR2743921}. 
(Note that $\mathbf{v}_{i, 0}$ in Fig.~1 of \cite{PannetonLM06} corresponds to this variable $\mathbf{v}_i$.) 
This approach is a key technique for drastically improving $N_1$ and $\Delta$. 
By refining the recursion formulas proposed in \cite{MR2743921}, 
we implement the state transition $f$ with the following recursions: 
\begin{eqnarray}
\mathbf{v}_{i+1} & := & (\overline{\mathbf{w}}_{i}^{w-r} \mid \underline{\mathbf{w}}_{i+1}^{r} ) A \oplus \mathbf{w}_{i+M} \oplus \mathbf{v}_{i} B \label{eqn:MELG_recurtion1_def},\\
\mathbf{w}_{i+N-1} & : = &  (\overline{\mathbf{w}}_{i}^{w-r} \mid \underline{\mathbf{w}}_{i+1}^{r} ) \oplus \mathbf{v}_{i+1} C. \label{eqn:MELG_recurtion2_def}
\end{eqnarray}
$A$, $B$, and $C$ are $(w \times w)$-matrices defined indirectly as follows:
\begin{eqnarray}
\mathbf{w}A & := & \begin{cases} 
(\mathbf{w} \gg1) & \text{if $w_{w-1} = 0$,} \label{eqn:def_mat_A}\\
(\mathbf{w} \gg1) \oplus \mathbf{a} & \text{if $w_{w-1} = 1$,}
\end{cases}\\
\mathbf{w}B & := & \mathbf{w} \oplus (\mathbf{w} \ll \sigma_1), \label{eqn:def_mat_B} \\
\mathbf{w}C & := & \mathbf{w} \oplus (\mathbf{w} \gg \sigma_2) \label{eqn:def_mat_C},
\end{eqnarray}
where $\mathbf{w} = (w_0, \ldots, w_{w-1}) \in \mathbb{F}_2^w$ and $\mathbf{a} \in \mathbb{F}_2^w$ are $w$-bit vectors, 
$\sigma_1$ and $\sigma_2$ are integers with $0 < \sigma_1, \sigma_2< w$, and 
``$\mathbf{w} \ll l$'' and  ``$\mathbf{w} \gg l$'' denote left and right logical (i.e., zero-padded) shifts by $l$ bits, respectively. 
Note that $\mathbf{v}_i$ is one component in a state $s_i \in S$ and is not the output. 

\begin{figure}
\centering
\includegraphics[width=9cm]{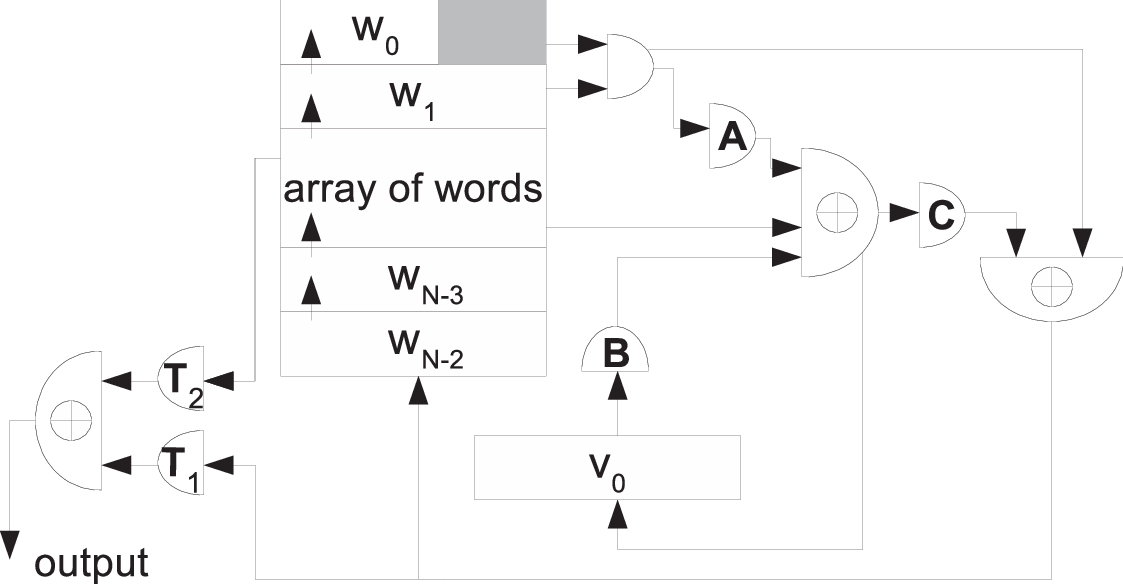}
  \caption{Circuit-like description of the proposed generators.}
  \label{fig:MELG}
\end{figure}

To attain the maximal equidistribution exactly, 
we design a linear output transformation using another word in the state array, which comes from \cite{Harase2009}. 
More precisely, for the right-hand side of (\ref{eqn:MELG_state_transition}), 
we consider the following linear output transformation $o: S \to O$ with one more memory reference: 
\begin{eqnarray} \label{eqn:new_tempering}
(\overline{\mathbf{w}}_{i+1}^{r-w}, \mathbf{w}_{i +2}, \ldots, \mathbf{w}_{i + N-1}, \mathbf{v}_{i+1}) \in S \mapsto \mathbf{w}_{i+N-1} T_1 \oplus \mathbf{w}_{i+L} T_2 \in O. 
\end{eqnarray}
Here $T_1$ and $T_2$ are $(w \times w)$-matrices defined by  
\begin{eqnarray}
\mathbf{w}T_1 & := & \mathbf{w} \oplus (\mathbf{w} \ll \sigma_3), \label{eqn:def_mat_T1}\\
\mathbf{w}T_2 & := & (\mathbf{w} \ \& \ \mathbf{b}), \label{eqn:def_mat_T2}
\end{eqnarray}
where $L$ is an integer with $0 < L < N-2$, $\sigma_3$ is an integer with $0 < \sigma_3 < w$, 
$\&$ denotes bitwise AND, and $\mathbf{b} \in \mathbb{F}_2^w$ is a $w$-bit vector. 
A circuit-like description of the proposed generators is shown in Figure~\ref{fig:MELG}. 

An equivalent formal algorithm can be described as Algorithm~1, in which, 
instead of shifting, we use a round-robin technique (i.e., a pointer technique) 
to improve the efficiency of the generation. 
Let $\mathbf{w}[0..N-2]$ be an array of $N-1$ unsigned integers with $w$ bits,
and $\mathbf{v}$ be a $w$-bit unsigned integer that corresponds to the extra state variable $\mathbf{v}_i$. 
Let $\mathbf{x}$ be a temporary variable 
and $\mathbf{y}$ be an output variable that are $w$-bit unsigned integers, respectively. 
Set $\mathbf{m}^{w-r} \gets (\underset{w-r}{\underbrace{1, \ldots, 1}}, \underset{r}{\underbrace{0, \ldots, 0}})$ and 
$\mathbf{\hat{m}}^{r} \gets (\underset{w-r}{\underbrace{0, \ldots, 0}}, \underset{r}{\underbrace{1, \ldots, 1}})$. 
Here $\mathbf{m}^{w-r}$ is a bit mask that retains the first $w -r$ bits
and sets the other $r$ bits to zero, whereas $\mathbf{\hat{m}}^{r}$ is its bitwise complement. Before Algorithm~1 begins, we initialize $(\mathbf{w}[0] \ \& \ \mathbf{m}^{w-r}), \mathbf{w}[1], \ldots, \mathbf{w}[N-2], \mathbf{v} \gets \mbox{initial values, not all zero}$, and set the pointer $i \gets 0$. 
For more detailed descriptions of the initialization, see Remark~\ref{rmk:initialization}. 

\begin{algorithm}[htbp] \label{algo:MELG}
\SetAlgoNoLine
\DontPrintSemicolon

\caption{The algorithm of the proposed generators}
 \SetKwInOut{Input}{Input}
 \SetKwInOut{Output}{Output}
 \Input{$\mathbf{w}[0], \mathbf{w}[1], \ldots, \mathbf{w}[N-2], \mathbf{v}$, and the pointer $i$}
\Output{$w$-bit output $\mathbf{y}$}
Set $\mathbf{x} \gets (\mathbf{w}[i] \ \& \ \mathbf{m}^{w-r}) \oplus (\mathbf{w}[(i +1) \mod{(N-1)}] \ \& \ \mathbf{\hat{m}}^{r}).$ // compute $(\overline{\mathbf{w}}_{i}^{w-r} \mid \underline{\mathbf{w}}_{i+1}^{r} )$.\;
Set $\mathbf{v} \gets \mathbf{x}A \oplus \mathbf{w}[(i + M) \mod{(N-1)}] \oplus \mathbf{v}B.$ // compute Eq.~(\ref{eqn:MELG_recurtion1_def}).\;
Set $\mathbf{w}[i] \gets \mathbf{x} \oplus \mathbf{v}C$. // compute Eq.~(\ref{eqn:MELG_recurtion2_def}).\;
Set $\mathbf{y} \gets \mathbf{w}[i]T_1 \oplus \mathbf{w}[(i+L) \mod{(N-1)}]T_2$. // compute Eq.~(\ref{eqn:new_tempering}).\;
Increment $i \gets i + 1$. If $i \geq N-1$, then $i \gets i \mod{(N-1)}$. // increment the pointer $i$.\;Return $\mathbf{y}$\;
\label{alg:one}
\end{algorithm}

\subsection{Specific Parameters} \label{subsec:parameters}
We search for specific parameters in the following way. 
First, we look for $M$, $\sigma_1$, $\sigma_2$, and $\mathbf{a} \in \mathbb{F}_2^w$ in (\ref{eqn:MELG_recurtion1_def})--(\ref{eqn:def_mat_C}) at random
such that $f$ attains the maximal period $2^p-1$. 
In general, because we can obtain several parameters, 
we choose a parameter set 
whose $N_1$ is large enough and whose output has large $k(v)$ for $v = 1, \ldots, w$ as far as possible 
in the case where we set $\mathbf{y} \gets \mathbf{w}[i]$ in Algorithm~1 (i.e.,  $T_1$ and $T_2$ are the identity and zero matrices, respectively). 
In the next step, we search for $L$, $\sigma_3$, and $\mathbf{b} \in \mathbb{F}_2^w$ in (\ref{eqn:new_tempering})--(\ref{eqn:def_mat_T2}) at random 
such that the generator is ``almost" maximally equidistributed (i.e., $\Delta$ is almost $0$). 
Finally, to obtain $\Delta = 0$ strictly, we apply a slight modification to the bit mask $\mathbf{b}\in \mathbb{F}_2^{w}$ by using the backtracking algorithm in \cite{Harase2009} (with some trial and error). 
Table~\ref{table:MELG} lists specific parameters for  64-bit maximally equidistributed 
generators with periods ranging from $2^{607}-1$ to $2^{44497}-1$. 
We attach the acronym 64-bit MELGs, for 64-bit ``maximally equidistributed $\mathbb{F}_2$-linear generators" with Mersenne prime period, to the proposed generators. 
The code in C is available at \url{https://github.com/sharase/melg-64}.

\begin{table}[t]
\begin{center}
\tbl{Specific Parameters of 64-bit Maximally Equidistributed $\mathbb{F}_2$-Linear Generators (MELGs) \label{table:MELG}}{
\begin{tabular}{cccccc} \hline
{} & $M$ & $\sigma_1$ & $\sigma_2$ & $\mathbf{a}$ & $N_1$\\
{} & $L$ & $\sigma_3$ & {} & $\mathbf{b}$ & $\Delta$ \\ \hline \hline
\multicolumn{6}{l}{$p = 607, w = 64, r = 33, N = 10$} \\ \hline
MELG607-64 & $5$ & $13$ & $35$ & ${\tt 81f1fd68012348bc}$ & $313$ \\
{} & $3$ & $30$ & {} & ${\tt 66edc62a6bf8c826}$ & $0$ \\ \hline \hline
\multicolumn{6}{l}{$p = 1279, w = 64, r = 1, N = 20$} \\ \hline
MELG1279-64 & $7$ & $22$ & $37$ & ${\tt 1afefd1526d3952b}$ & $641$\\
{} & $5$ & $6$ & {} & ${\tt 3a23d78e8fb5e349}$ & $0$ \\ \hline \hline
\multicolumn{6}{l}{$p = 2281, w = 64, r = 23, N = 36$} \\ \hline
MELG2281-64 & $17$ & $36$ & $21$ & ${\tt 7cbe23ebca8a6d36}$ & $1145$ \\
{} & $6$ & $6$ & {} & ${\tt e4e2242b6e15aebe }$ & $0$ \\ \hline \hline
\multicolumn{6}{l}{$p = 4253, w = 64, r = 35, N = 67$} \\ \hline
MELG4253-64 & $29$ & $30$ & $20$ & ${\tt fac1e8c56471d722}$ & $2129$ \\
{} & $9$ & $5$ & {} & ${\tt cb67b0c18fe14f4d}$ & $0$ \\ \hline \hline
\multicolumn{6}{l}{$p = 11213, w = 64, r = 51, N = 176$} \\ \hline
MELG11213-64 & $45$ & $33$ & $13$ & ${\tt ddbcd6e525e1c757}$ & $5455$ \\
{} & $4$ & $5$ & {} & ${\tt bd2d1251e589593f}$ & $0$ \\ \hline \hline
\multicolumn{6}{l}{$p = 19937, w = 64, r = 31, N = 312$} \\ \hline
MELG19937-64 & $81$ & $23$ & $33$ & ${\tt 5c32e06df730fc42}$ & $9603$ \\
{} & $19$ & $16$ & {} & ${\tt 6aede6fd97b338ec}$ & $0$ \\ \hline \hline
\multicolumn{6}{l}{$p = 44497, w = 64, r = 47, N = 696$} \\ \hline
MELG44497-64 & $373$ & $37$ & $14$ & ${\tt 4fa9ca36f293c9a9}$ & $19475$ \\
{} & $95$ & $6$ & {} & ${\tt 6fbbee29aaefd91}$ & $0$ \\ \hline \hline
\end{tabular}}
\end{center}
\end{table}

In parallel computing, an important requirement is the availability of pseudorandom number generators with disjoint streams. 
These are usually implemented by partitioning 
the output sequences of a long-period generator into long disjoint subsequences whose starting points are found by making large jumps in the original sequences.  
For this purpose, \citeN{MR2437205} proposed a fast jumping-ahead algorithm for $\mathbb{F}_2$-linear generators. 
We also implemented this algorithm for our $64$-bit MELGs. 
The code is also available at the above website. The default skip size is $2^{256}$. 

For our generators, we implemented a function to produce double-precision floating-point numbers $u_0, u_1, u_2, \ldots$ in $[0,1)$ in IEEE 754 format 
by using the method of Section~2 of \cite{MR2743921} (i.e., the 12 bits for sign and exponent are kept constant, and the 52 bits of the significand are taken from the generator output). We note that this method is preferable from the viewpoint of $k(v)$ because it does not introduce approximation errors by division. 

\begin{remark} \label{rmk:initialization}
\citeN{Matsumoto:2007:CDI:1276927.1276928} reported that many pseudorandom number generators have some nonrandom bit patterns when initial seeds are systematically chosen, especially when the initialization scheme is based on a linear congruential generator. To avoid such phenomena, the 2002 version of MT19937 adopted a nonlinear initializer described in Eq.\ (30) of \cite{Matsumoto:2007:CDI:1276927.1276928}. For our proposed generators, 
we implemented a similar initialization scheme.  
Let $\mathbf{w}_0 \in \mathbb{F}_2^{w}$ be an initial seed. 
To obtain an initial state $s_0 = (\overline{\mathbf{w}}_0^{w-r}, \mathbf{w}_{1}, \mathbf{w}_{2}, \ldots, \mathbf{w}_{N-2}, \mathbf{v}_0) \in S$, 
we set $\overline{\mathbf{w}}_0^{w-r}$ as the $w-r$ most significant bits of $\mathbf{w}_0$ 
and compute
\[ \mathbf{w}_{i} \gets a \times (\mathbf{w}_{i-1} \oplus (\mathbf{w}_{i-1} \gg \sigma_4)) + i \pmod{2^w} \] 
for $i = 1, \ldots, N-2$ and 
\[ \mathbf{v}_0 \gets a \times (\mathbf{w}_{N-2} \oplus (\mathbf{w}_{N-2} \gg \sigma_4)) + N-1 \pmod{2^w}, \] 
using the usual integer arithmetic $+, -$, and $\times$, 
where $a = 6364136223846793005$ (in decimal notation) is a multiplier recommended in \cite[pp.~104]{Knuth:1997:ACP:270146}, $\sigma_4 = 62$, and $w= 64$. 
These parameters are from the 2004 version of 64-bit MT \cite{Nishimura2000} mentioned in Section~\ref{sec:performance}. 
In a similar manner, we implemented an initializer whose seed is an array of integers of arbitrary length. 
\end{remark}

\begin{remark}
As pointed out in Section~3.8 of \cite{MR3031631}, it might be difficult to run the proposed generators 
in a multi-threaded environment, such as GPUs. This is because our generators have the heavy dependencies on the partial computation of the extra state variable $\mathbf{v}_i$ in the recursion (\ref{eqn:MELG_recurtion1_def}). 
Thus, in this case, it seems that the original MTs \cite{MT19937,Nishimura2000} are suitable 
because of the simplicity of recursion. We note that MTGP is designed specifically for this purpose (see \cite{MR3031631} for details). 
\end{remark}

\section{Performance} \label{sec:performance}

We compare the following $\mathbb{F}_2$-linear generators corresponding to 64-bit integer output sequences:
\begin{itemize}
\item MELG19937-64: the 64-bit integer output of our proposed generator;
\item MT19937-64: the 64-bit integer output of the 64-bit Mersenne Twister (downloaded from \url{http://www.math.sci.hiroshima-u.ac.jp/~m-mat/MT/emt64.html}); 
\item MT19937-64 (ID3): the 64-bit integer output of a 64-bit Mersenne Twister based on a five-term recursion (ID3) \cite{Nishimura2000};
\item SFMT19937-64 (without SIMD): the 64-bit integer output of the SIMD-oriented Fast Mersenne Twister SFMT19937 without SIMD \cite{SFMT};
\item SFMT19937-64 (with SIMD): the 64-bit integer output of the foregoing with SIMD \cite{SFMT}.
\end{itemize}
The first three generators have period length $2^{19937}-1$. 
SFMT19937-64 has the period of a multiple of $2^{19937}-1$ (see Proposition~1 in \cite{SFMT} for details). 
Table~\ref{table:benchmark} summarizes the figures of merit $N_1$, $\Delta$, and timings. 
In this table, we report the CPU time (in seconds) taken to generate $10^9$ 64-bit unsigned integers for each generator. 
The timings were obtained using two 64-bit CPUs: (i) a 3.40 GHz Intel Core i7-3770 and (ii) a 2.70 GHz Phenom II X6 1045T.
The code was written in C and compiled with GCC using the -O3 optimization flag on 64-bit Linux operating systems. 
For SFMT19937-64, we measured the CPU time for the case of sequential generation (see \cite{SFMT} for details). 

\begin{table}[t]
\tbl{Figures of Merit $N_1, \Delta$, and CPU Time (in Seconds) Taken to Generate $10^9$ 64-bit Unsigned Integers \label{table:benchmark}}{
  \centering
  \begin{tabular}{|l|c|c|c|c|} \hline 
    Generators & $N_1$ & $\Delta$ & CPU time (Intel) & CPU time (AMD) \\ \hline \hline
    MELG19937-64 & $9603$ & $0$ & $4.2123$ & $6.2920$\\
    MT19937-64 & $285$ & $7820$ & $5.1002$ & $6.6490$\\ 
    MT19937-64 (ID3) & $5795$ & $7940$ & $4.8993$ & $6.7930$\\ 
    SFMT19937-64 (without SIMD) & $6711$ & $14095$ & $4.2654$ & $5.6123$\\ \hline \hline
    SFMT19937-64 (with SIMD) & $6711$ & $14095$ & $1.8457$ & $2.8806$ \\ \hline
  \end{tabular}}
\end{table} 

MELG19937-64 is maximally equidistributed 
and also has $N_1 \approx p/2$. These values are the best in this table. 
In terms of generation speed, MELG19937-64 is comparable to or even slightly faster than the MT19937-64 generators on the above two platforms.  
SFMT19937-64 without SIMD is comparable to or faster than MELG19937-64, and
SFMT19937-64 with SIMD is more than twice as fast as MELG19937-64. 
However, $\Delta$ for SFMT19937-64 is rather large. 
In fact, the SFMT generators are optimized under the assumption that one will mainly be using 32-bit output sequences, 
so that the dimensions of equidistribution with $v$-bit accuracy for 64-bit output sequences are worse than those for 32-bit cases ($\Delta = 4188$). 
For this, we analyze the structure of SFMT19937 in the online Appendix. 

Finally, we convert 64-bit integers into double-precision floating-point numbers $u_0, u_1, u_2, \ldots$ in $[0,1)$,
and submit them to the statistical tests included in the SmallCrush, Crush, and BigCrush batteries of TestU01~\cite{MR2404400}. 
Note that these batteries have 32-bit resolution and have not yet been tailored to 64-bit integers. 
In our C-implementation, for MELG19937-64 and MT19937-64, 
we generate the above uniform real numbers by {\tt (x >> 11) * (1.0/9007199254740992.0)}, where {\tt x} is 
a 64-bit unsigned integer output. 
For SFMT19937-64, we use a function {\tt sfmt\_genrand\_res53()}, which is obtained by dividing 64-bit unsigned integers ${\tt x}$ by $2^{64}$, 
i.e., {\tt x * (1.0/18446744073709551616.0)}; 
see the online Appendix for details. 
In any case, we investigate the 32 most significant bits of 64-bit outputs in TestU01. 
The generators in Table~\ref{table:MELG} passed all the tests; except for the linear complexity tests (unconditional failure) 
and matrix-rank tests (failure only for small $p$), 
which measure the $\mathbb{F}_2$-linear dependency of the outputs and reject $\mathbb{F}_2$-linear generators. 
This is a limitation of $\mathbb{F}_2$-linear generators. 
However, for some matrix-rank tests, we can observe differences between MELGs and SFMTs. 
Table~\ref{table:MatrixRank} summarizes the $p$-values on the matrix-rank test of No.~60 of Crush for five initial states.  
(SFMT1279-64, SFMT2281-64, and SFMT19937-64 denote the results for double-precision floating-point numbers in $[0,1)$ converted from the 64-bit integer outputs of SFMT1279, SFMT2281, and SFMT19937 in \cite{SFMT}, respectively.)

\begin{table}[t]
\centering
\tbl{$p$-Values on the Matrix-Rank Test of No.~60 of Crush in TestU01 \label{table:MatrixRank}}{
\begin{tabular}{l||c|c|c|c|c} \hline
{ } & {1st} & {2nd} & {3rd} & {4th} & {5th} \\ \hline  \hline
MELG1279-64 & $0.89$ & $0.41$ & $0.11$ & $0.70$ & $0.22$ \\ \hline 
MELG2281-64 & $0.70$ & $0.02$ & $0.62$ & $0.49$ & $0.98$ \\ \hline 
MELG19937-64 & $0.23$ & $0.13$ & $0.32$ & $0.14$ & $0.85$ \\ \hline 
SFMT1279-64 & $<10^{-300}$ & $<10^{-300}$ & $< 10^{-300}$ &  $<10^{-300}$ & $< 10^{-300}$ \\ \hline 
SFMT2281-64 & $<10^{-300}$ & $<10^{-300}$ & $<10^{-300}$ & $<10^{-300}$ & $<10^{-300}$ \\ \hline 
SFMT19937-64 & $0.29$ & $0.02$ & $0.06$ & $0.49$ & $0.83$ \\ \hline 
\end{tabular}}
\end{table}

\begin{remark} \label{rmk:reverse}
We occasionally see the use of the least significant bits of pseudorandom numbers in applications. 
An example is the case in which one generates uniform integers from $0$ to $15$ 
by taking the bit mask of the 4 least significant bits or modulo $16$. 
For this, we invert the order of the bits (i.e., the $i$-th bit is exchanged with the $(w-i)$-th bit) in each integer 
and compute the dimension of equidistribution with $v$-bit accuracy, dimension defect at $v$, and total dimension defect for inversion, 
which are denoted by $k'(v)$, $d'(v)$, and $\Delta'$, respectively. 
In this case, MELG19937-64 is not maximally equidistributed, 
but $\Delta' = 4047$ and $d'(v)$ is $0$ or $1$ for each $v \leq 11$. 
Note that $\Delta'$ takes values $9022$, $8984$, $21341$ for MT19937-64, MT19937-64 (ID3), SFMT19937-64, respectively.  
$\Delta'$ of MELG19937-64 is still smaller than $\Delta$ of the other generators in Table~\ref{table:benchmark}. 
However, as far as possible, we recommend using the most significant bits (e.g., by taking the right-shift in the above example), 
because our generators are optimized preferentially from the most significant bits.
\end{remark}
\begin{remark} \label{rmk:conversion} 
For 32-bit generators, there have been some implementations that produce 64-bit unsigned integers or 53-bit double-precision floating-point numbers 
(in IEEE 754 format)
by concatenating two consecutive 32-bit unsigned integers. 
We note that such conversions might not be preferable from the viewpoint of $k(v)$. 
As an example, consider the 32-bit MT19937 generator in the header {\tt <random>} of the C++11 STL in GCC (see \cite[Chapter~26.5]{ISO:2012:III}). 
Let $\mathbf{z}_0, \mathbf{z}_1, \mathbf{z}_2, \ldots \in \mathbb{F}_2^{32}$ 
be a 32-bit unsigned integer sequence from 32-bit MT19937. 
To obtain 64-bit unsigned integers, the GCC implements a random engine adaptor {\tt independent\_bit\_engine} to produce 64-bit unsigned integers from 
the concatenations as
\begin{eqnarray} \label{eqn:concatenation1}
(\mathbf{z}_{0}, \mathbf{z}_{1}), (\mathbf{z}_{2}, \mathbf{z}_{3}), (\mathbf{z}_{4}, \mathbf{z}_{5}),  (\mathbf{z}_{6}, \mathbf{z}_{7}), \ldots \in \mathbb{F}_2^{64}. 
\end{eqnarray}
To generate 53-bit double-precision floating-point numbers in $[0, 1)$ (i.e., {\tt uniform\_real\_distribution(0,1)} for MT19937), 
the GCC implementation generates 64-bit unsigned integers 
\begin{eqnarray} \label{eqn:concatenation2}
(\mathbf{z}_{1}, \mathbf{z}_{0}), (\mathbf{z}_{3}, \mathbf{z}_{2}), (\mathbf{z}_{5}, \mathbf{z}_{4}),  (\mathbf{z}_{7}, \mathbf{z}_{6}), \ldots \in \mathbb{F}_2^{64} 
\end{eqnarray}
by concatenating two consecutive 32-bit integer outputs and divides them by the maximum value $2^{64}$. 
The sequences (\ref{eqn:concatenation1}) and (\ref{eqn:concatenation2}) can be viewed as $\mathbb{F}_2$-linear generators 
with the state transition $f^2$ in (\ref{eqn:MT_state_transition}), 
so that we can compute $k(v)$. 
Note that the sequences (\ref{eqn:concatenation1}) and (\ref{eqn:concatenation2}) are different: (\ref{eqn:concatenation2}) is obtained by exchanging the 32 most significant bits for the 32 least significant bits in each 64-bit word in (\ref{eqn:concatenation1}), that is, the $\mathbb{F}_2$-linear output functions are described as $s_i \in S \mapsto (o(s_i), o(f(s_i))) \in \mathbb{F}_2^{64}$ in (\ref{eqn:concatenation1}) and $s_i \in S \mapsto (o(f(s_i)), o(s_i)) \in \mathbb{F}_2^{64}$ in (\ref{eqn:concatenation2}). In fact, their $k(v)$'s are different for $v = 33, \ldots, 64$. As a result, we have $\Delta = 13543$ for $v = 1, \ldots, 64$ in (\ref{eqn:concatenation1}) and 
$\Delta = 13161$ for $v = 1, \ldots, 52$ in (\ref{eqn:concatenation2}), which are worse than $\Delta$ of MT19937-64. 
In particular, $k(12) = 623 < \lfloor 19937/12 \rfloor  = 1661$ for each case. 
For this reason, we feel that there is a need to design 64-bit high-quality pseudorandom number generators. 
\end{remark}
\section{Conclusions} \label{sec:conclusions}

In this article, we have designed 64-bit maximally equidistributed $\mathbb{F}_2$-linear generators with Mersenne prime period   
and searched for specific parameters with period lengths from $2^{607}-1$ to $2^{44497}-1$. 
The key techniques are (i) state transitions with double feedbacks and (ii) linear output transformations with several memory references. 
As a result, the generation speed is still competitive with 64-bit Mersenne Twisters on some platforms. 
The code in C is available at \url{https://github.com/sharase/melg-64}. 
Pseudorandom number generation is a trade-off between speed and quality. 
Our generators offer both high performance and computational efficiency. 



\begin{acks}
The authors would like to thank Associate Editor Professor Pierre L'Ecuyer and the anonymous reviewers 
for their valuable comments and suggestions. The authors also wish to
express their gratitude to Professor Makoto Matsumoto at Hiroshima University for many helpful comments 
and Professor Syoiti Ninomiya at Tokyo Institute for continuous encouragement. 
\end{acks}
\bibliographystyle{ACM-Reference-Format-Journals}
\bibliography{harase-kimoto-bib}

\received{May 2015}{May 2016, March 2017}{November 2017}


\end{document}